\documentclass[11pt]{article}
\usepackage{mathrsfs}
\usepackage{latexsym,amsmath,amssymb,amsfonts,epsfig}
\usepackage{color,colordvi,amscd,indentfirst,graphics}
\allowdisplaybreaks

\newtheorem{theorem}{Theorem}[section]

\newtheorem{lemma}[theorem]{Lemma}
\newtheorem{corollary}[theorem]{Corollary}

\begin{document}

\title{A note on independence number, connectivity and $k$-ended tree}
\author{Pham Hoang Ha\footnote{E-mail address: ha.ph@hnue.edu.vn.}\\
Department of Mathematics\\
Hanoi National University of Education\\
136 XuanThuy Street, Hanoi, Vietnam\\}
\date{}

\maketitle{}

\bigskip

\begin{abstract}
A $k$-ended tree is a tree with at most $k$ leaves. In this note, we give a simple proof for the following theorem. Let $G$ be a connected graph and $k$ be an integer ($k\geq 2$). Let $S$ be a vertex subset of $G$ such that $\alpha_{G}(S) \leq  k + \kappa_{G}(S)- 1.$ Then, $G$ has a $k$-ended tree which covers $S.$ Moreover, the condition is sharp.
\end{abstract}




\section{Introduction} 
In this note, we only consider finite simple graphs. Let $G$ be a
graph with vertex set $V(G)$ and edge set $E(G)$. A subset $X\subseteq V(G)$ is called an \emph{independent set} of
$G$ if no two vertices of $X$ are adjacent in $G$.  The maximum size
of independent sets in $G$ is denoted by $\alpha(G)$. A graph G is $k$-connected if it has more than $k$ vertices
and every subgraph obtained by deleting fewer than $k$ vertices is connected; the connectivity of $G,$
written $\kappa(G),$ is the maximum $k$ such that $G$ is $k$-connected.  For any $S\subseteq V(G)$, we
denote by $|S|$ the cardinality of $S$. We define
$\alpha_G(S)$ the maximum cardinality of independent sets of $S$ in $G,$ which is called the
independence number of $S$ in $G.$ For two vertices $x, y$ of $G,$ the local connectivity
$\kappa_G(x, y)$ is defined to be the maximum number of internally disjoint paths connecting
$x$ and $y$ in $G.$ We define $\kappa_G(S) := \min\{\kappa_G(x, y) : x, y \in S, x \not= y\}.$ Moreover, if
$| S| = 1,$ $\kappa_G(S)$ is defined to be $+\infty.$ When $S=G,$ we have $\alpha_G(G)=\alpha(G)$  and by Menger's theorem we have $\kappa_G(S)=\kappa(G).$ A Hamiltonian cycle (path) is a cycle (path) which passes through all vertices of a
graph.

In 1972, Chv\'{a}tal and Erd\H{o}s proved the following famous theorem which related to the independence number, connectivity and Hamiltonian cycle (path) of a graph.
\begin{theorem}[{\cite[Chv\'{a}tal and Erd\H{o}s]{CE}}] Let G be a connected graph.
	\item [{(1)}] If $\alpha (G) \leq \kappa(G),$ then G has a Hamiltonian cycle unless $G = K_1$ or $K_2$.
	\item [{(2)}] If $\alpha (G) \leq \kappa(G)+1,$ then G has a Hamiltonian path.
\end{theorem}

Let $T$ be a tree. A vertex of degree one is a \emph{leaf} of $T$
and a vertex of degree at least three is a \emph{branch vertex} of
$T$. A tree having at most $k$ leaves is called $k$-ended tree. Then a Hamiltonian path is nothing but a spanning $2$-ended tree. In 1979, Win improved the above result by proving the following theorem.
\begin{theorem}[{\cite[Win]{Wi79}}]\label{theo1.1}
	Let $G$ be a graph and let $k$ be an integer ($k\geq 2$).  If
	$\alpha(G)\leq k+\kappa(G)-1$, then $G$ has a spanning tree with at most $k$ leaves.
\end{theorem}

On the other hand, when we consider a cycle (path) containing specified vertices of a graph as a
generalization of a Hamiltonian cycle (path), many results were invented.
\begin{theorem}[{\cite[Fournier]{F} }] Let $G$ be a $2$-connected graph, and let $S \subseteq V(G).$ If
	$\alpha_G (S) \leq \kappa(G),$ then G has a cycle covering $S.$
\end{theorem}
\begin{theorem} [{\cite[Ozeki and Yamashita]{OY}}] Let $G$ be a $2$-connected graph and let
	$S \subseteq V(G).$ If $\alpha_G (S) \leq \kappa_G(S)$, then $G$ has a cycle covering $S.$
\end{theorem}

A natural question is whether Win's result can be improved by giving a sharp condition to show the existence of the $k$-ended tree covering a given subset of $V(G).$ In this note, we give an affirmative answer to this question. In particular, we prove the following theorem. 
\begin{theorem}\label{main}
	Let $G$ be a connected graph and $k$ be an integer ($k\geq 2$). Let $S$ be a subset of $V(G)$ such that $\alpha_{G}(S) \leq  k + \kappa_{G}(S)- 1.$ Then, $G$ has a $k$-ended tree covering $S.$ 
\end{theorem}
It is easy to see that if a tree has at most $k$ leaves ($k\geq 2$),
then it has at most $k-2$ branch vertices. Therefore, we immediately
obtain the following corollary from Theorem~\ref{main}.
\begin{corollary}\label{coro}
	Let $G$ be a connected graph and $k$ be an integer ($k\geq 2$). Let $S$ be a subset of $V(G)$ such that $\alpha_{G}(S) \leq  k + \kappa_{G}(S)- 1.$ Then, $G$ has a tree $T$ such that $T$ covers $S$ and has at most $k-2$ branch vertices. 
\end{corollary}

We first show that the conditions of Theorem \ref{main} and Corollary \ref{coro} are sharp. Let  $m,k\geqslant 1$ be integers, and let $K_{m,m+k}=(A, B)$ be a complete bipartite graph with $|A|=m, |B|=m+k$. Set $S= B.$ Then we are easy to see that $\alpha_{G}(S) =  k + \kappa_{G}(S)$  and every tree covering $S$ has at most $k+1$ leaves. Moreover it also has at most $k-1$ branch vertices. Therefore, the conditions of Theorem \ref{thm-main} and Corollary \ref{coro} are sharp.

To prove Theorem \ref{main}, we prove a slightly stronger following result.
\begin{theorem}\label{thm-main}
	Let $G$ be a connected graph and $k$ be an integer ($k\geq 2$). Let $S$ be a subset of $V(G).$ Then either $G$ has a $k$-ended tree $T$ covering $S,$ or there exists a $k$-ended tree $T$ in $G$ such that 
	\begin{equation*}
		\alpha_{G}(S-V(T)) \leq \alpha_{G}(S)  - \kappa_{G}(S)- k+ 1.
	\end{equation*} 
\end{theorem}

 Beside that many researches on the relations of independence number, connectivity and the tree whose maximum degree is at most $k$ containing specified
vertices of a graph are studied. We would like to refer the readers the papers \cite{CMOT}, \cite{OY},\cite{OY11} and \cite{Yan} for more details.
\section{Proof of Theorem \ref{thm-main}}
By using the same technique in \cite{K}, Yan in \cite{Yan} proved the following result. It needs for the proof of Theorem \ref{thm-main}.
\begin{lemma}[{\cite[Corollary 1]{Yan}}]\label{lem1}
	Let $G$ be a connected graph and $S \subseteq V(G).$ Then either the vertices of $S$ can be covered by one path of $G,$ or there exists a path $P$ of $G$ such that 
	\begin{equation*}
		\alpha_{G}(S-V(P)) \leq \alpha_{G}(S) - \kappa_{G}(S)- 1.
	\end{equation*} 
\end{lemma}
Next, we prove Theorem \ref{thm-main} by induction on $k(\geq 2).$

For $k=2,$ by Lemma \ref{lem1}, the theorem holds.

Assume that the theorem  holds for some $k=t\geq 2,$ that is, either the vertices of $S$ can be covered by one $t$-ended tree of $G,$ or there exists a $t$-ended tree $T$ of $G$ such that
\begin{equation}\label{eq1}
	\alpha_{G}(S-V(T)) \leq \alpha_{G}(S) -\kappa_{G}(S)- t+ 1.
\end{equation} 

If there exists a $(t+1)$-ended tree such that it covers $S$ then the theorem holds for $k=t+1.$ Otherwise, every $(t+1)$-ended tree of $G$ does not cover $S.$ In particular, $S$ can not be covered by any $t$-ended tree of $G.$ By the induction hypothesis, there exists a $t$-ended tree $T$ of $G$ such that (\ref{eq1}) is correct. Let $S_1, ..., S_m$ be all subsets of $S-V(T)$ such that $|S_i|=\alpha_{G}(S-V(T))$ for all $i\in\{1,...,m\}.$ For each vertex $s\in \cup_{i=1}^mS_i,$ since $G$ is connected, there exists some path joining $s$ to $T.$ Denote by $P[s, T]$ the set of such paths in $G.$ We choose a maximal path $P_{0}$ in $\{P[s, T]| s \in \cup_{i=1}^{m} S_i \}.$ Assume that $P_{0}$ joins the vertex $s_0\in \cup_{i=1}^mS_i$ to $T.$ Now, we prove that $P_{0}\cap S_i \not= \emptyset$ for all $i\in \{1,...,m\}.$ Indeed, otherwise, there exists some $j$ such that $P_{0}\cap S_j = \emptyset.$ By $|S_j|=\alpha_{G}(S-V(T))$ and $s_0 \in S-V(T),$ there exists some vertex $s_j\in S_j$ such that $s_{0}s_{j}\in E(G).$ We consider the path $P'=P_0+s_{0}s_{j}.$ Then $P'$ joins $s_j$ to $T$ and $|P'| > |P_0|,$ which implies a contradiction with the maximality of $P_0.$ Therefore we conclude that  $P_{0}\cap S_i \not= \emptyset$ for all $i\in \{1,...,m\}.$ Now, we set $T'=T+P_{0}.$ Then $T'$ has at most $(t+1)$ leaves. On the other hand, because $P_{0}\cap S_i \not= \emptyset$ and $|S_i|=\alpha_{G}(S-V(T))$ for all $i\in \{1,...,m\},$ we obtain $\alpha_{G}(S-V(T'))\leq \alpha_{G}(S-V(T))-1.$ So
 \begin{equation*}
	\alpha_{G}(S-V(T')) \leq \alpha_{G}(S-V(T)) -1\leq \alpha_{G}(S) -\kappa_{G}(S)- t.
\end{equation*}
This implies that the theorem holds for $k=t+1.$

Therefore, the theorem holds for all $k \geq 2$ by the principle of mathematical induction. Hence we complete the proof of Theorem \ref{thm-main}.
\bigskip

{\bf Acknowledgements.} The research is supported by the NAFOSTED Grant of Vietnam (No. 101.04-2018.03).

\end{document}